\documentclass[11pt]{amsart}
\usepackage[margin=1in]{geometry}
\usepackage[T1]{fontenc}
\usepackage{lmodern}
\usepackage{amsmath,amssymb}
\usepackage{microtype}
\usepackage[hidelinks]{hyperref}
\newtheorem{theorem}{Theorem}[section]
\newtheorem{lemma}[theorem]{Lemma}
\newtheorem{proposition}[theorem]{Proposition}
\newtheorem{corollary}[theorem]{Corollary}
\theoremstyle{remark}
\newtheorem*{example}{Example}
\numberwithin{equation}{section}
\newcommand{\simplex}{\overline{\Delta}^m}
\newcommand{\norm}[1]{\|#1\|_\infty}
\newcommand{\PP}{\mathbb P}
\newcommand{\E}{\mathbb E}
\title[Dice enterprises and uniform query tails]{Dice enterprises and polynomial approximation: uniform query tails}
\author{Jessica Taylor}
\email{jessica@mediangroup.org}
\date{October 2026}
\hypersetup{pdftitle={Dice enterprises and polynomial approximation: uniform query tails},
pdfauthor={Jessica Taylor}}

\begin{document}
\begin{abstract}
A binary-output dice enterprise uses independent rolls of a die with unknown face
probabilities $\mathbf p$ (along with fair coin flips) to simulate a coin with bias $f(\mathbf p)$.
We determine which power-law tails, uniform over the closed probability
simplex, are attainable for the number $N$ of rolls. For continuous $f$
with values in $(0,1)$ and $q>0$, some enterprise for $f$ satisfies
$\sup_{\mathbf p}\PP_{\mathbf p}(N>n)=O(n^{-q})$ if and only if $f$ can
be approximated uniformly by polynomials of degree $k$ with error
$O(k^{-2q})$. When this holds, the same tail bound is attained by an oblivious
enterprise, which chooses its number of rolls before rolling. The
characterization extends to multiple dice, with total degree on a
product of simplices; it then suffices to roll the dice in a fixed
cyclic order. For $q>1$, these enterprises have a uniformly bounded
expected number of rolls.
\end{abstract}
\maketitle

\section{Introduction and main results}

A \emph{Bernoulli factory} for a function $f$ turns independent tosses
of a coin with unknown bias $p$, together with auxiliary fair coin flips,
into a coin with bias $f(p)$. A \emph{dice enterprise} does the same
with a die whose face probabilities $\mathbf p$ are unknown. We study
the two-outcome case: given independent rolls of the die, produce a
Bernoulli random variable with parameter $f(\mathbf p)$. This is a special case of the
finite-output formulation of Morina et al.~\cite{morina2022dice}, and we
call such a procedure simply an \emph{enterprise} for $f$. Finite
outputs reduce to this case by a standard rejection argument
(Section~\ref{setup}), and a two-sided die recovers the classical
Bernoulli factory.

Beyond existence, it is natural to ask how many rolls an enterprise
needs. We seek bounds that hold uniformly over the unknown die weighting:
for which $q$ can an enterprise for $f$ achieve
$\sup_{\mathbf p}\PP_{\mathbf p}(N>n)=O(n^{-q})$? The answer is a
polynomial approximation rate for $f$, with the exponent doubled.

\subsection*{Setting}
For $m\geq0$, the parameter space of a die with $m+1$ faces is the
closed simplex
\[
\overline\Delta^m
=\Bigl\{\mathbf p=(p_0,\ldots,p_m)\in[0,1]^{m+1}:
\sum_{i=0}^m p_i=1\Bigr\}.
\]
For a compact set $S$ and an integer $k\geq0$, let $\mathcal P_k^S$ be
the space of restrictions to $S$ of real polynomials of total degree at
most $k$, and let
\[
E_k^S(f)=\inf_{P\in\mathcal P_k^S}\|f-P\|_{\infty,S}
\]
be the error of best uniform approximation. The degree of a restriction
is the least total degree of a polynomial representing it; on a simplex,
or on a product of simplices, this is its degree in affine coordinates.

An \emph{enterprise} for $f:\simplex\to[0,1]$ is a randomized procedure,
fixed in advance, that may request independent rolls of a die with weighting
$\mathbf p$ and independent fair coin flips, each request depending only on
what it has observed so far. For every $\mathbf p\in\simplex$ it halts
almost surely and outputs $Y\in\{0,1\}$ with
$\PP_{\mathbf p}(Y=1)=f(\mathbf p)$. Its \emph{query count} $N$ is the
number of rolls; fair coin flips are not counted. An enterprise is
\emph{oblivious} if it first uses fair coin flips alone to choose a finite
budget, and then rolls the die exactly that many times. The distribution of $N$
for an oblivious enterprise does not depend on $\mathbf p$, although the enterprise's
output may depend on all the rolls. The number of fair coin flips used must be
finite almost surely, but we bound neither its expectation nor the
computational cost of the enterprise; uncomputable decision trees are allowed. Section~\ref{setup}
gives a formal definition.

\emph{Conventions.} A uniform tail bound $O((1+n)^{-q})$ means that
$\PP_{\mathbf p}(N>n)\leq C(1+n)^{-q}$ for all
inputs $\mathbf p$ and integers $n\geq0$, with a
single constant $C$. We also write $O(n^{-q})$
when only the asymptotic rate is relevant. Likewise,
$E_k(f)=O(k^{-2q})$ means there is a single $C$ with $E_k(f) \leq Ck^{-2q}$ for all $k\geq1$.  Implicit constants, including those in $A\lesssim B$ (which
means $A\leq CB$ for nonnegative $A$ and $B$) and in $A\asymp B$ (which
means $A\lesssim B\lesssim A$), never depend on the input
probabilities or on indices that vary; they may depend on fixed data
such as the target $f$, the dimension $m$, the exponent $q$, and the enterprise at
hand. Subscripts on $\lesssim$ indicate dependence where it matters. For $f:[0,1]\to[0,1]$, an enterprise for $f$ means an enterprise
for $(1-p,p)\mapsto f(p)$, as in the usual Bernoulli factory convention.

\begin{theorem}[Uniform-tail characterization]\label{dice_approx_equiv}
Let $m\geq0$ and $q>0$, and let $f:\simplex\to(0,1)$ be continuous.
The following are equivalent:
\begin{enumerate}
\item[(i)] Some enterprise for $f$ satisfies
\[
\sup_{\mathbf p\in\simplex}\PP_{\mathbf p}(N>n)=O((1+n)^{-q}).
\]
\item[(ii)] The best approximation errors satisfy
\[
E_k^{\simplex}(f)=O(k^{-2q}).
\]
\end{enumerate}
When these hold, the tail bound in (i) is attained by an oblivious
enterprise. For $q>1$, this enterprise has
$\sup_{\mathbf p}\E_{\mathbf p}N<\infty$.
\end{theorem}

Adaptivity does not improve the set of attainable uniform power-tail
exponents. The domain includes the boundary of the simplex, where some faces have
probability zero. Since $f$ is continuous on a compact set and takes
values in $(0,1)$, there is $\epsilon\in(0,1/2)$ with
$f(\simplex)\subseteq[\epsilon,1-\epsilon]$; the construction uses this
margin. The implication (i)$\Rightarrow$(ii) holds for every
$[0,1]$-valued target, and in particular forces continuity
(Proposition~\ref{enterprise_to_approx}). The range restriction is needed
for the converse, even in the Bernoulli case: $f(p) = (p - 1/2)^2$ has no
enterprise. Indeed, at $p=1/2$ every leaf of the tree is reached with
positive probability, so $f(1/2)=0$ forces every leaf to be labelled $0$,
contradicting $f(0)=1/4$.
\begin{example}
Consider a Bernoulli input coin with unknown probability $p$ of
outputting $1$, and let $f(p)=(1+\sqrt p)/3$.
Every polynomial $P(p)$ of degree at most $k$ gives an even
polynomial $P(u^2)$ of degree at most $2k$, with the same
uniform approximation error for $|u|$ as $P$ has for
$\sqrt p$. Conversely, any degree-$2k$ approximant to $|u|$
can be made even by averaging its values at $u$ and $-u$,
without increasing the error, and every such even polynomial
has the form $P(u^2)$ with $\deg P\le k$. Thus Bernstein's
classical estimate~\cite{bernstein1914} gives
\[
E_k^{[0,1]}(f)
=\frac13 E_{2k}^{[-1,1]}(u\mapsto |u|)
\asymp k^{-1}.
\]
By Theorem~\ref{dice_approx_equiv}, $f$ has an enterprise with
uniform tail $O((1+n)^{-1/2})$, but none with a uniform tail
$O((1+n)^{-q})$ for any $q>1/2$. In particular, no enterprise
for $f$ has uniformly bounded expected query count: by
Markov's inequality, such a bound would imply a uniform
tail $O((1+n)^{-1})$.

For comparison, let $g(p)=(1+|p-\tfrac12|)/3$.
Unlike $f$, which is only $1/2$-H\"older, $g$ is Lipschitz,
though it has a kink at $p=1/2$.
Under the substitution $u=2p-1$, we have
$g(p)=1/3+|u|/6$, so Bernstein's classical estimate gives
\[
E_k^{[0,1]}(g)
=\frac16 E_k^{[-1,1]}(u\mapsto |u|)
\asymp k^{-1}.
\]
Thus, despite their different ordinary smoothness,
both functions have the same largest attainable uniform
power-tail exponent, $q=1/2$: an interior Lipschitz kink
has the same effect on this exponent as an endpoint
square-root cusp.
\end{example}

\subsection*{Multiple dice}
Let $\ell\geq1$ and $m_1,\ldots,m_\ell\geq0$, and put
\[
K=\prod_{i=1}^{\ell}\overline\Delta^{m_i}.
\]
A multi-die enterprise may choose, at each step and depending on the
history so far, which die to roll; a query to die $i$ returns an
independent sample with weighting $\mathbf p_i$. Now $N$ is the total number
of queries to all dice. A multi-die enterprise is oblivious if, before
any query, it uses fair coin flips to choose a complete finite schedule of die
indices, which it then follows exactly. The number and dimensions of the
dice are fixed, and constants may depend on them.

\begin{theorem}[Multiple input dice]\label{multi_dice}
Let $q>0$ and let $f:K\to(0,1)$ be continuous. The following are
equivalent:
\begin{enumerate}
\item[(i)] Some multi-die enterprise for $f$ satisfies
\[
\sup_{(\mathbf p_1,\ldots,\mathbf p_\ell)\in K}
\PP_{\mathbf p_1,\ldots,\mathbf p_\ell}(N>n)=O((1+n)^{-q}).
\]
\item[(ii)] With degree meaning total degree in all coordinates of the
product,
\[
E_k^K(f)=O(k^{-2q}).
\]
\end{enumerate}
When these hold, the tail bound in (i) is attained by an oblivious
enterprise that rolls the dice in the cyclic order
$1,\ldots,\ell,1,\ldots,\ell,\ldots$; in particular, it rolls every die equally many times. For
$q>1$, its expected total query count is uniformly bounded on $K$.
\end{theorem}

The proof uses a joint-outcome die and its marginal map. This preserves the attainable tail exponents while controlling polynomial degrees, even when the original enterprise chooses its input die adaptively.

\subsection*{Relation to earlier work}
Keane and O'Brien~\cite{keane1994bernoulli} characterized the functions
that admit a Bernoulli factory. Nacu and
Peres~\cite[Theorem 2 and Proposition 3]{nacu2005fast} related
exponentially decaying tails to analyticity and developed a formulation
of simulation in terms of polynomial approximation. Holtz, Nazarov, and
Peres~\cite[Theorem 8]{holtz2011newcoins} showed that, for
$f:[0,1]\to(0,1)$ and noninteger $s>0$, $f\in C^s[0,1]$ if and only
if there is a factory and a constant $C$ such that
\[
\PP_p(N>n)\leq C\Delta_n(p)^s
\quad\text{for all }p\in[0,1]\text{ and }n\geq1,
\qquad
\Delta_n(p)=\max\left \{\sqrt{p(1-p)/n},1/n \right\}.
\]
Their envelope depends on the bias $p$: it is of order
$n^{-s/2}$ at fixed interior points and $n^{-s}$ at the endpoints. In contrast, Theorem~\ref{dice_approx_equiv} characterizes bounds
of the form $O(n^{-q})$ that hold uniformly over all input
probabilities, without requiring faster decay near the endpoints. Holtz et al. count rounds containing one input-coin toss and one fair-coin toss, whereas our query count charges only input rolls; thus these bounds use different conventions for auxiliary randomness.

For example, the function $f(p)=(1+\sqrt p)/3$ belongs
to $C^{1/2}[0,1]$. Applying the theorem of Holtz et al.\
with $s=1/2$ guarantees a tail of order $n^{-1/4}$
at fixed interior points, whereas our characterization
gives a uniform $O(n^{-1/2})$ tail.

Morina, \L atuszy\'nski, Nayar, and Wendland~\cite{morina2022dice} constructed dice enterprises for rational maps between open simplices, using the input die as the source of randomness. They show how to extract auxiliary fair coin flips from such a die. Here we instead permit independent fair coin flips without charging them as input queries, allowing a uniform query-tail analysis on the closed simplex, including its vertices.

Paes Leme and
Schneider~\cite{paesleme2023multiparameter} characterized
realizability for multiple unknown coins on the closed hypercube.
Nonadaptive factories have also been studied; see Koskela et
al.~\cite{koskela2026debiased} and the earlier constructions discussed
there.

These prior results do not directly characterize uniform power
tails through unconstrained best polynomial approximation.
The compatible Bernstein-polynomial approximations of
Nacu and Peres impose conditions beyond a bound on
$E_k(f)$, while the bias-dependent envelope of Holtz,
Nazarov, and Peres requires faster decay near the endpoints
than a uniform power-tail bound.

Our contribution is the exact equivalence between uniform
$O((1+n)^{-q})$ query tails and the approximation rate
$E_k(f)=O(k^{-2q})$, on closed simplices and their
products. Bernstein coefficient inversion gives controlled
coefficients for degree-$k$ corrections at representation
degree $O(k^2)$, permitting their simulation through a
mixture construction. Conversely, positive Bernstein
compression bounds the best approximation errors of
termination polynomials in terms of their uniform norms
and representation degrees. Together these estimates
establish the equivalence, with an oblivious enterprise
attaining every admissible exponent.

\subsection*{Outline.} Section \ref{construction} constructs oblivious enterprises from polynomial approximations, and Section 4 proves the converse using positivity of termination polynomials. Section 5 extends the characterization to multiple input dice, and Section 6 gives smoothness consequences.

\section{The model and Bernstein representations}\label{setup}

If $m=0$, the simplex is a point and the target $f$ is constant, which a
fixed-bias coin realizes with no queries; all approximation errors
vanish, and Theorem~\ref{dice_approx_equiv} is immediate. We therefore
assume $m\geq1$ throughout Sections~\ref{setup}--\ref{necessity}.

\subsection*{Enterprises as trees}
Formally, an enterprise is a rooted tree, fixed independently of
$\mathbf p$, whose internal vertices are input queries, with one child
per face, or fair-coin queries, with two children, and whose leaves are
labelled $0$ or $1$; its vertices are finite histories. Run on input
$\mathbf p$, it moves from an input-query vertex to the child of face
$i$ with probability $p_i$, and from a fair-coin vertex to either child
with probability $1/2$. For every $\mathbf p$, a leaf must be reached
after finitely many steps almost surely. The query count $N$ is the
number of input-query vertices visited, and $Y$ is the label of the
leaf reached.

We freely use coins with fixed biases and draws from fixed countable
distributions. Read the fair coin flips as the binary expansion of a
uniform point of $[0,1]$, and locate it in a partition into intervals
of the desired probabilities. Almost surely, a finite prefix determines
the interval. The endpoints may be noncomputable; the model does not
require computability.

\subsection*{Multiple outcomes}
For existence and uniform power-tail bounds, finite-output enterprises
reduce to the binary-output case.
Morina~\cite[Proposition~3.21]{morina2021extending} gives this reduction
with preservation of exponential tails (``fast simulation'').

Suppose $F_1,\ldots,F_d$ ($d \geq 2$) sum to 1 and each has a
binary enterprise. In each independent round, choose $i$ uniformly from
$\{1,\ldots,d\}$, run the enterprise for $F_i$,
and output $i$ if it returns $1$.  Each round succeeds
with probability $1/d$, so we have $\PP_{\mathbf{p}}(\mathrm{output}~i) = \frac{F_i(\mathbf{p})/d}{1/d} = F_i(\mathbf{p})$ as desired. 
Conversely, a finite-output enterprise gives a binary
enterprise for each $F_i$ by reporting whether its
output is $i$, with the same query count.

Suppose each binary enterprise has uniform tail
$O((1+n)^{-q})$, where $q>0$. Generate independent rounds, with $T_j$ the query count of round $j$, and
$S_j\in\{0,1\}$ the success indicator of round $j$. Let $R = \min\{j ~ | ~ S_j=1\}$ be the number
of rounds performed. Conditional on
$R=r$, the first $r-1$ rounds failed and the last
succeeded. By independence of rounds, for each $j\le r$,
\[
\PP_{\mathbf p}(T_j>t\mid R=r)
=
\PP_{\mathbf p}(T_j>t\mid S_j=s_j),
\]
where $s_j=0$ for $j<r$ and $s_r=1$. Since success and
failure each have probability at least $1/d$,
\[
\PP_{\mathbf p}(T_j>t\mid R=r)
\le d\,\PP_{\mathbf p}(T_j>t)
\lesssim (1+t)^{-q}.
\]
Thus, since $N=\sum_{j=1}^R T_j$, the union bound gives
\[
\PP_{\mathbf p}(N>n\mid R=r)
\lesssim r(1+n/r)^{-q}
\le r^{q+1}(1+n)^{-q}.
\]
Averaging over $R$ yields
\[
\PP_{\mathbf p}(N>n)
\lesssim \mathbb E[R^{q+1}](1+n)^{-q}
=O((1+n)^{-q}),
\]
since $R$ is geometric and hence has finite moments
of every order. This reduction need not preserve obliviousness.

\subsection*{Bernstein representations}
For $n\geq0$, let $\Lambda_n^m=\{\mathbf a\in\mathbb N_0^{m+1}:|\mathbf a|=n\}$,
and write $\mathbf p^{\mathbf a}=\prod_i p_i^{a_i}$, with $0^0=1$, and
$\binom n{\mathbf a}=n!/\prod_i a_i!$. Every $P\in\mathcal P_n^{\simplex}$
has a unique \emph{degree-$n$ Bernstein representation}
\begin{equation}\label{bernstein_rep}
P(\mathbf p)=\sum_{\mathbf a\in\Lambda_n^m}
c(\mathbf a)\binom n{\mathbf a}\mathbf p^{\mathbf a},
\end{equation}
with \emph{Bernstein coefficients} $c(\mathbf a)$. The
\emph{representation degree} $n$ may exceed the degree of $P$. For
$I\subseteq\mathbb R$, let $\mathcal P_{n,I}^{\simplex}$ be the set of
polynomials whose degree-$n$ Bernstein coefficients all lie in $I$.

The basis polynomials in \eqref{bernstein_rep} are nonnegative and sum
to one, so $\norm P\leq\max_{\mathbf a}|c(\mathbf a)|$. If all
coefficients lie in $[0,1]$, then $P$ has an oblivious enterprise with
exactly $n$ queries: roll the die $n$ times, let $\mathbf a$ be the
vector of face counts, and output $1$ with probability $c(\mathbf a)$.
This works because the count vector equals $\mathbf a$ with probability
$\binom n{\mathbf a}\mathbf p^{\mathbf a}$. All our constructions are
mixtures of such enterprises.

For $n\geq1$, the \emph{Bernstein operator} is
\begin{equation}\label{bernstein_operator}
B_nh(\mathbf p)=\E_{\mathbf p}h(\hat{\mathbf p}_n)
=\sum_{\mathbf a\in\Lambda_n^m}
h(\mathbf a/n)\binom n{\mathbf a}\mathbf p^{\mathbf a},
\end{equation}
where $\hat{\mathbf p}_n$ is the empirical distribution of $n$ rolls.
It is a contraction in the uniform norm. Unless another domain is
indicated, norms and approximation errors refer to $\simplex$.

\subsection*{Weighted edge derivatives}
Let $e_0,\ldots,e_m$ be the standard basis of $\mathbb R^{m+1}$. The
unit edge directions of the simplex, and their weights, are
\[
V_{\simplex}=\{(e_i-e_j)/\sqrt2:0\leq i<j\leq m\},
\qquad
\phi_{(e_i-e_j)/\sqrt2}(\mathbf p)=\sqrt{2p_ip_j}.
\]
For $\xi\in V_{\simplex}$, $\partial_\xi$ is the derivative along $\xi$,
which is tangent to the simplex, and $\phi_\xi^r\partial_\xi^rh$ denotes
the weight $\phi_\xi^r$ times the derivative $\partial_\xi^rh$. Writing
$\partial_i-\partial_j$ for the derivative along $e_i-e_j$, we have
$\phi_\xi^2\partial_\xi^2=p_ip_j(\partial_i-\partial_j)^2$ for
$\xi=(e_i-e_j)/\sqrt2$, so the operator
\begin{equation}\label{diffusion_operator}
\mathcal LP=\sum_{\xi\in V_{\simplex}}\phi_\xi^2\partial_\xi^2P
=\sum_{0\leq i<j\leq m}p_ip_j(\partial_i-\partial_j)^2P
\end{equation}
is the operator in Voronovskaya's formula $n(B_nh-h)\to\mathcal Lh/2$
for smooth $h$, reflecting that $\hat{\mathbf p}_n-\mathbf p$ has
covariance matrix $n^{-1}(p_i\delta_{ij}-p_ip_j)_{i,j}$.
Appendix~\ref{coordinates} identifies these weighted
derivatives with those of Ditzian and Zhou~\cite{ditzian1993optimal};
the identification preserves degree and uniform norm on the full closed
simplex. For $C^r(\simplex)$ we take restrictions of $C^r$ functions
defined on a neighborhood of $\simplex$ in its affine hull; we apply
these notions only to polynomials.

\section{From approximation to an oblivious enterprise}\label{construction}

A polynomial with Bernstein coefficients in $[0,1]$ is realized by the
elementary enterprise of Section~\ref{setup}. A good approximant $P$ of
degree $k$ need not have such coefficients, but its degree-$n$
Bernstein coefficients approach the grid values $P(\mathbf a/n)$,
relative to $\norm P$, as $n/k^2\to\infty$. The mechanism is as
follows. If $B_n$ is
invertible on polynomials of degree at most $k$, then
$P=B_n(B_n^{-1}P)$, and by \eqref{bernstein_operator} the degree-$n$
Bernstein coefficients of $P$ are precisely the values of
$B_n^{-1}P$ on the grid $\{\mathbf a/n\}$. When $k^2\ll n$, the
operator $B_n$ is close to the identity on such polynomials, and so
is its inverse. Thus these coefficients are close to
$P(\mathbf a/n)$.

\begin{lemma}[Bernstein inversion]\label{bern_promote}
There is a constant $C_m\geq2$, depending only on $m$, such that for all
integers $k,n\geq1$ and all $P\in\mathcal P_k^{\simplex}$,
\begin{equation}\label{bernstein_poly}
\norm{B_nP-P}\leq C_m\frac{k^2}{n}\norm P.
\end{equation}
Moreover, $B_n$ maps $\mathcal P_k^{\simplex}$ into itself. If
$\eta=C_mk^2/n<1$, then the degree-$n$ Bernstein coefficients of $P$
satisfy
\begin{equation}\label{coefficient_bound}
\max_{\mathbf a\in\Lambda_n^m}|c(\mathbf a)-P(\mathbf a/n)|
\leq\frac{\eta}{1-\eta}\norm P.
\end{equation}
In particular, if $\eta\leq1/2$, then
$\max_{\mathbf a}|c(\mathbf a)|\leq2\norm P$.
\end{lemma}

\begin{proof}
Suppose first that $k^2<n$. Use the change of coordinates in Appendix~\ref{coordinates}. Ditzian and Zhou's Voronovskaya-type
estimate~\cite[Theorem~5.1, (5.1)]{ditzian1993optimal}
gives
\[
\Bigl\|B_nP-P-\frac1{2n}\mathcal LP\Bigr\|_\infty
\lesssim_m\frac{k^4}{n^2}\norm P.
\]
Their weighted Bernstein
inequality~\cite[(3.1)]{ditzian1993optimal},
applied with derivative order $2$ and
$p=\infty$, gives
\[
\norm{\mathcal LP}\lesssim_m k^2\norm P.
\]

 Since
$k^4/n^2<k^2/n$, this proves \eqref{bernstein_poly} when $k^2<n$. When
$k^2\geq n$, the uniform-norm contraction property gives
$\norm{B_nP-P}\leq2\norm P\leq2(k^2/n)\norm P$.

To see that $B_n$ preserves $\mathcal P_k^{\simplex}$, let $X_1,\ldots,X_n$
be the rolls, write
$\hat p_{n,i}=n^{-1}\sum_{t=1}^n\mathbf 1_{\{X_t=i\}}$, and expand a
monomial of degree $d\leq k$ in these empirical frequencies. Each term
is a product of $d$ indicators. It vanishes if some roll is required to
show two different faces; otherwise it involves at most $d$ distinct
rolls, so by independence its expectation is a monomial in $\mathbf p$
of degree at most $d$.

Thus $B_n$ is an operator on the finite-dimensional normed space
$\mathcal P_k^{\simplex}$, equipped with the uniform norm.
The bound \eqref{bernstein_poly} gives
$\|I-B_n\|_{\mathrm{op}}\leq\eta$, where
$\|\cdot\|_{\mathrm{op}}$ is the induced operator norm.
If $\eta<1$, the Neumann series
$B_n^{-1}=\sum_{j\geq0}(I-B_n)^j$ converges in operator norm and gives
\[
\|B_n^{-1}-I\|_{\mathrm{op}}\leq\frac{\eta}{1-\eta},
\qquad
\|B_n^{-1}\|_{\mathrm{op}}\leq\frac{1}{1-\eta}.
\]
Since $\eta<1$ and $C_m\geq2$, we have $n>k$, so the
degree-$n$ Bernstein representation is available. Put
$\widetilde P=B_n^{-1}P\in\mathcal P_k^{\simplex}$. Since
$P=B_n\widetilde P$, formula \eqref{bernstein_operator} writes $P$ in
the degree-$n$ Bernstein basis, with coefficients
$c(\mathbf a)=\widetilde P(\mathbf a/n)$. Therefore
$|c(\mathbf a)-P(\mathbf a/n)|\leq\norm{\widetilde P-P}
\leq\eta\norm P/(1-\eta)$, which is \eqref{coefficient_bound}, and
$|c(\mathbf a)|\leq\norm{\widetilde P}\leq\norm P/(1-\eta)$.
\end{proof}

Raising the representation degree makes the coefficients of a small
correction small, but they may have either sign. A baseline polynomial
whose coefficients stay away from $0$ and $1$ leaves room for
corrections of both signs.

\begin{lemma}[Mixing signed corrections]\label{correct_obliv}
Let $n_0,n_1,\ldots$ be nonnegative integers, and let
$\gamma_1,\gamma_2,\ldots\geq0$ satisfy $A=\sum_{j\geq1}\gamma_j<1/2$.
Suppose $R\in\mathcal P_{n_0,[A,1-A]}^{\simplex}$ and
$Q_j\in\mathcal P_{n_j,[-\gamma_j,\gamma_j]}^{\simplex}$ for $j\geq1$.
Then $R+\sum_{j\geq1}Q_j$ has an oblivious enterprise that first
chooses branch $0$ with probability $1-2A$ and branch $j$ with
probability $2\gamma_j$ for each $j\geq1$. On branch $j\geq0$ it uses
exactly $n_j$ queries.
\end{lemma}

\begin{proof}
Since $\norm{Q_j}\leq\gamma_j$, the series converges absolutely and
uniformly. Discard the indices with $\gamma_j=0$, for which $Q_j=0$, and
put
\[
W_0=\frac{R-A}{1-2A},\qquad
W_j=\frac12+\frac{Q_j}{2\gamma_j}\quad(j\geq1).
\]
Each $W_j$ has degree-$n_j$ Bernstein coefficients in $[0,1]$, so it is
realized with exactly $n_j$ queries. Choose a branch $j\geq0$ with the
stated probabilities, using fair coin flips only, and run the enterprise for
$W_j$. The branch probabilities sum to one, and the output is $1$ with
probability
\[
(1-2A)W_0+\sum_{j\geq1}2\gamma_jW_j
=R-A+\sum_{j\geq1}(\gamma_j+Q_j)=R+\sum_{j\geq1}Q_j.
\qedhere
\]
\end{proof}

\begin{proposition}[Approximation gives an oblivious tail]\label{approx_to_enterprise}
Let $q>0$ and $\epsilon\in(0,1/2)$, and let
$f:\simplex\to[\epsilon,1-\epsilon]$ be continuous with
$E_k(f)=O(k^{-2q})$. Then $f$ has an oblivious enterprise with
$\sup_{\mathbf p}\PP_{\mathbf p}(N>n)=O((1+n)^{-q})$.
\end{proposition}

\begin{proof}
\emph{Approximants.} Choose $C>0$ with $E_k(f)\leq Ck^{-2q}$ for all
$k\geq1$, and put $\alpha_j=C4^{-qj}$ for $j\geq0$. Best approximants
exist in the finite-dimensional spaces $\mathcal P_k^{\simplex}$, so
there are $R_j\in\mathcal P_{2^j}^{\simplex}$ with
$\norm{f-R_j}\leq\alpha_j$. For $j\geq1$, put
\[
Q_j=R_j-R_{j-1},\qquad
\gamma_j=2(\alpha_j+\alpha_{j-1}),
\]
so that $\norm{Q_j}\leq\gamma_j/2$ and $\sum_j\gamma_j<\infty$. Fix $J$
so large that
\[
\alpha_J\leq\epsilon/4
\qquad\text{and}\qquad
A:=\sum_{j>J}\gamma_j<\epsilon/2.
\]
Since $R_j\to f$ uniformly, $f=R_J+\sum_{j>J}Q_j$.

\emph{Baseline.} The values of $R_J$ lie in $[3\epsilon/4,1-3\epsilon/4]$.
By Lemma~\ref{bern_promote}, applied to the fixed polynomial $R_J$, we
can choose a representation degree $n_{\mathrm{base}}$ for which the
right side of \eqref{coefficient_bound} is at most $\epsilon/4$. The
coefficients then lie in
$[\epsilon/2,1-\epsilon/2]\subseteq[A,1-A]$.

\emph{Corrections.} Fix an integer $L\geq2C_m$ and put $n_j=L4^j$ for
$j>J$. Since $Q_j$ has degree at most $2^j$, Lemma~\ref{bern_promote}
applies with $\eta=C_m4^j/n_j\leq1/2$, and bounds the degree-$n_j$
coefficients of $Q_j$ in absolute value by $2\norm{Q_j}\leq\gamma_j$.

\emph{The enterprise.} Lemma~\ref{correct_obliv}, with indices shifted,
now gives an oblivious enterprise for $f$. The baseline branch has
probability $1-2A$ and budget $n_{\mathrm{base}}$; correction branch
$j>J$ has probability $2\gamma_j$ and budget $L4^j$. Thus correction $j$ pairs
degree $2^j$ with query budget of order $4^j$ and probability of order
$4^{-qj}$. For $n\geq\max\{n_{\mathrm{base}},L4^{J+1}\}$, choose
$j\geq J+1$ with $L4^j\leq n<L4^{j+1}$. Then
\[
\PP_{\mathbf p}(N>n)=2\sum_{i>j}\gamma_i
\lesssim\sum_{i\geq j}4^{-qi}
\lesssim4^{-qj}\lesssim n^{-q},
\]
using $\gamma_i\lesssim4^{-q(i-1)}$ and $n<L4^{j+1}$. Enlarging the
constant covers the finitely many smaller $n$. The bound is uniform in $\mathbf p$ because the distribution
of the query count $N$ does not depend on $\mathbf p$.
\end{proof}

\section{From query tails to approximation}\label{necessity}

The probability $\PP_{\mathbf p}(Y=1,\ N\leq n)$ is a polynomial
of degree at most $n$ (Lemma~\ref{segment_bernstein}) and differs from
$f$ by at most $\PP_{\mathbf p}(N>n)$. For a tail $O(n^{-q})$
this only shows $E_n(f)=O(n^{-q})$, half the desired exponent. The
missing information is positivity: the probability of
outputting $1$ within a block of query counts has nonnegative Bernstein
coefficients (Lemma~\ref{segment_bernstein}), and such polynomials can
be approximated at degree of order $\sqrt n$ to any fixed relative accuracy (Lemma~\ref{approx_pos}).

\begin{lemma}[Positive termination blocks]\label{segment_bernstein}
For any enterprise and integers $0\leq n'\leq n$, the function
\[
\mathbf p\longmapsto\PP_{\mathbf p}(Y=1,\ n'\leq N\leq n)
\]
belongs to $\mathcal P_{n,[0,1]}^{\simplex}$.
\end{lemma}

\begin{proof}
Sample a string $\mathbf X=(X_1,\ldots,X_n)$ of rolls in advance, and
feed $X_1,X_2,\ldots$ to the enterprise as it requests rolls. Runs that
request more than $n$ rolls do not contribute to the event, so the event
is determined by $\mathbf X$ and the fair coin flips. For a fixed string
$\mathbf x$, let $b(\mathbf x)\in[0,1]$ be the probability, over the
fair coin flips alone, that the enterprise run on $\mathbf x$ halts with
output $1$ after between $n'$ and $n$ rolls. Since the tree does not
depend on $\mathbf p$, neither does $b(\mathbf x)$; runs that never halt
on a particular fixed string simply do not count. A string with count
vector $\#\mathbf x=\mathbf a$ has probability $\mathbf p^{\mathbf a}$,
so grouping strings by count vector gives the representation
\eqref{bernstein_rep}, with
\[
c(\mathbf a)=\binom n{\mathbf a}^{-1}
\sum_{\substack{\mathbf x\in\{0,\ldots,m\}^n\\
\#\mathbf x=\mathbf a}}b(\mathbf x)\in[0,1].
\]
This holds also when some face probabilities vanish.
\end{proof}

\subsection{Two weighted estimates and positive compression}

We use two standard estimates. The first is Lorentz's inequality~\cite{lorentz1963degree} for
polynomials with nonnegative Bernstein coefficients: for every integer
$r\geq1$, a univariate polynomial $H$ with nonnegative degree-$n$
Bernstein coefficients, $n\geq1$, satisfies
\begin{equation}\label{lorentz}
\sup_{0\leq u\leq1}
\bigl|(u(1-u))^{r/2}H^{(r)}(u)\bigr|
\lesssim_r n^{r/2}\|H\|_{\infty,[0,1]}.
\end{equation}
We derive this from the form given by Holtz, Nazarov, and
Peres~\cite[Result 11, (8)]{holtz2011newcoins}. If
$M=\|H\|_{\infty,[0,1]}>0$, take $\Omega(h)=M(1+h)$ there. It is
positive and increasing, satisfies $\Omega(2h)\leq2\Omega(h)$, and
satisfies $|H(u)|\leq\Omega(\Delta_n(u))$ for
$\Delta_n(u)=\max\{\sqrt{u(1-u)/n},1/n\}$. Since
$\Delta_n\leq1$, the cited inequality gives
$|H^{(r)}(u)|\leq2\kappa_rM\Delta_n(u)^{-r}$ for a constant $\kappa_r$,
and multiplying by $(u(1-u))^{r/2}\leq n^{r/2}\Delta_n(u)^r$ gives
\eqref{lorentz}. The case $M=0$ is trivial.

The second is a weighted Jackson estimate: for integers $r,k\geq1$ and
$h\in C^r(\simplex)$,
\begin{equation}\label{jackson}
E_k(h)\lesssim_{m,r}k^{-r}
\Bigl(\norm h+
\max_{\xi\in V_{\simplex}}\norm{\phi_\xi^r\partial_\xi^r h}\Bigr).
\end{equation}
This follows from the weighted Jackson theorem, in the form stated by
Ditzian and Zhou~\cite[Theorem 3.3, (3.7)--(3.8)]{ditzian1993optimal},
which, in the
coordinates of Appendix~\ref{coordinates}, gives
$E_k(h)\lesssim_{m,r}K_r(h,k^{-r})+k^{-r}\norm h$ for the $K$-functional
\[
K_r(h,t^r)=\inf_{g\in C^r(\simplex)}
\Bigl(\norm{h-g}+t^r
\max_{\xi\in V_{\simplex}}\norm{\phi_\xi^r\partial_\xi^r g}\Bigr).
\]
Taking $g=h$, so that $K_r(h,k^{-r})\le k^{-r}\max_{\xi\in V_{\simplex}}
\norm{\phi_\xi^r\partial_\xi^r h}$, gives \eqref{jackson}.

\begin{lemma}[Positive Bernstein compression]\label{approx_pos}
For every integer $r\geq1$ there is $C_{m,r}>0$, depending only on $m$
and $r$, such that for all integers $n,k\geq1$ and all
$P\in\mathcal P_{n,[0,\infty)}^{\simplex}$,
\begin{equation}\label{compression}
E_k(P)\leq C_{m,r}\norm P
\min\Bigl\{1,\Bigl(\frac{\sqrt n}{k}\Bigr)^r\Bigr\}.
\end{equation}
\end{lemma}

\begin{proof}
We bound the weighted derivatives of $P$ by restricting to edge
segments. Fix $\mathbf p\in\simplex$ and $\xi=(e_i-e_j)/\sqrt2\in V_{\simplex}$,
and put $s=p_i+p_j$. If $s=0$, then $\phi_\xi(\mathbf p)=0$ and the
weighted derivative vanishes at $\mathbf p$. Otherwise, redistribute
the mass $s$ between faces $i$ and $j$: for $0\leq u\leq1$, let
$\zeta(u)\in\simplex$ have coordinates $\zeta(u)_i=su$ and
$\zeta(u)_j=s(1-u)$, with all other coordinates equal to those of
$\mathbf p$, and put $H(u)=P(\zeta(u))$. Then
$\|H\|_{\infty,[0,1]}\leq\norm P$.

The polynomial $H$ has nonnegative degree-$n$ Bernstein coefficients.
Indeed, each term $c(\mathbf a)\binom n{\mathbf a}\mathbf p^{\mathbf a}$
of the representation of $P$ restricts to a nonnegative multiple of
$u^{a_i}(1-u)^{a_j}$, and multiplying by $(u+(1-u))^{n-a_i-a_j}=1$
writes it in degree $n$ with nonnegative coefficients. Thus
\eqref{lorentz} applies to $H$.

Since $\zeta'(u)=s\sqrt2\,\xi$ and
$\phi_\xi(\zeta(u))=s\sqrt2\,(u(1-u))^{1/2}$, the chain rule gives
\begin{equation}\label{slice_identity}
\phi_\xi(\zeta(u))^r\,\partial_\xi^rP(\zeta(u))
=(u(1-u))^{r/2}(s\sqrt2)^r\partial_\xi^rP(\zeta(u))
=(u(1-u))^{r/2}H^{(r)}(u);
\end{equation}
the factors depending on the length of the segment cancel. At
$u=p_i/s$ we have $\zeta(u)=\mathbf p$, so \eqref{lorentz} yields
$|\phi_\xi(\mathbf p)^r\partial_\xi^rP(\mathbf p)|\lesssim_r n^{r/2}\norm P$.
The bound is uniform over $\mathbf p$ and $\xi$, including points where
$p_i$ or $p_j$ vanishes, so
\[
\max_{\xi\in V_{\simplex}}\norm{\phi_\xi^r\partial_\xi^rP}
\lesssim_r n^{r/2}\norm P.
\]
By \eqref{jackson} and $n\geq1$, it follows that
$E_k(P)\lesssim_{m,r}k^{-r}n^{r/2}\norm P$. Together with the trivial
bound $E_k(P)\leq\norm P$, this proves \eqref{compression}.
\end{proof}

The representation degree $n$ in Lemma~\ref{approx_pos} may exceed
the polynomial degree of $P$. Moreover, the bound scales with
$\norm P$, which can be much smaller than the largest Bernstein
coefficient. Thus termination blocks of small uniform probability
have correspondingly small approximation errors.
\subsection{Summing termination blocks}

\begin{proposition}[Query tails give approximation]\label{enterprise_to_approx}
Let $q>0$, and suppose an enterprise for $f:\simplex\to[0,1]$ satisfies
$\sup_{\mathbf p}\PP_{\mathbf p}(N>n)=O((1+n)^{-q})$. Then $f$ is
continuous and $E_k(f)=O(k^{-2q})$.
\end{proposition}

\begin{proof}
Split the event of output $1$ into blocks of query counts at the scales
$4^j$:
\[
\Pi_0(\mathbf p)=\PP_{\mathbf p}(Y=1,\ N\leq1),\qquad
\Pi_j(\mathbf p)=\PP_{\mathbf p}(Y=1,\ 4^{j-1}<N\leq4^j)\quad(j\geq1).
\]
By Lemma~\ref{segment_bernstein}, $\Pi_j\in\mathcal P_{4^j,[0,1]}^{\simplex}$
for every $j\geq0$. For $j\geq1$, the tail bound gives
\[
\norm{\Pi_j}\leq C(1+4^{j-1})^{-q}\leq C4^q4^{-qj},
\]
and $\norm{\Pi_0}\leq1$; hence $\norm{\Pi_j}\lesssim4^{-qj}$ for all
$j\geq0$. Since the enterprise halts almost surely,
$f=\sum_{j\geq0}\Pi_j$, and the series converges uniformly. In
particular, $f$ is continuous.

Fix an integer $r>2q$. Lemma~\ref{approx_pos}, with representation
degree $4^j$, gives
\[
E_k(\Pi_j)\lesssim 4^{-qj}\min\{1,(2^j/k)^r\}.
\]
Given $k\geq1$, let $J=\lfloor\log_2 k\rfloor$. Approximate each $\Pi_j$
with $j\leq J$ by a best approximant of degree $k$, and the remaining
blocks by $0$. The sum of these approximants has degree at most $k$, so
\[
E_k(f)\lesssim k^{-r}\sum_{j=0}^{J}2^{(r-2q)j}
+\sum_{j>J}2^{-2qj}\lesssim k^{-2q}.
\]
Here $2^J\leq k<2^{J+1}$. Since $r>2q$, the first sum is dominated by
its last term, of order $k^{r-2q}$, and the second is a geometric tail
of order $k^{-2q}$. With $r=2q$, the first sum would have $J+1$
comparable terms, costing a logarithmic factor.
\end{proof}

\begin{proof}[Proof of Theorem~\ref{dice_approx_equiv}]
The case $m=0$ was settled at the start of Section~\ref{setup}. For
$m\geq1$, the implication (i)$\Rightarrow$(ii) is
Proposition~\ref{enterprise_to_approx}, which applies to arbitrary
adaptive enterprises. For (ii)$\Rightarrow$(i), compactness provides
the range margin $\epsilon$, and Proposition~\ref{approx_to_enterprise}
provides an oblivious enterprise with the required tail. When $q>1$,
summing the tail bound gives
\[
\sup_{\mathbf p}\E_{\mathbf p}N
=\sup_{\mathbf p}\sum_{n=0}^{\infty}\PP_{\mathbf p}(N>n)
\leq C\sum_{n=0}^{\infty}(1+n)^{-q}<\infty.
\qedhere
\]
\end{proof}

\section{Multiple input dice}\label{multiple_inputs}

We deduce Theorem~\ref{multi_dice} from Theorem~\ref{dice_approx_equiv}
by replacing the $\ell$ dice with a single die whose outcomes are
tuples. The reduction must control both the approximation degree and
adaptive choices of which die to roll; it is not enough to view product
distributions as a subset of the joint simplex.

Let $\mathcal D=\prod_{i=1}^{\ell}\{0,\ldots,m_i\}$, and let
$S=\overline\Delta^{|\mathcal D|-1}$ be the simplex of all probability
distributions on $\mathcal D$, not necessarily of product form. The
marginal map $\mu:S\to K$ and the product map $\tau:K\to S$ are
\begin{equation}\label{product_maps}
\mu(\mathbf z)_{i,j}=\sum_{\mathbf a\in\mathcal D:\,a_i=j}z_{\mathbf a},
\qquad
\tau(\mathbf p_1,\ldots,\mathbf p_\ell)_{\mathbf a}
=\prod_{i=1}^{\ell}p_{i,a_i}.
\end{equation}
The map $\mu$ is linear, $\tau$ is polynomial of degree at most $\ell$,
and $\mu\circ\tau=\mathrm{id}_K$; in particular, $\mu$ is onto. Let
$g=f\circ\mu$, a continuous function on all of $S$ with values in
$(0,1)$.

\begin{lemma}[Approximation sandwich]\label{product_approx}
For every integer $k\geq0$,
\begin{equation}\label{product_sandwich}
E_{\ell k}^K(f)\leq E_k^S(g)\leq E_k^K(f).
\end{equation}
Consequently, for every $\beta>0$,
$E_k^K(f)=O(k^{-\beta})$ if and only if $E_k^S(g)=O(k^{-\beta})$.
\end{lemma}

\begin{proof}
If $P$ has degree at most $k$ on $K$, then $P\circ\mu$ has degree at
most $k$ on $S$ and, since $\mu$ is onto, the same uniform error:
$\|g-P\circ\mu\|_{\infty,S}=\|f-P\|_{\infty,K}$. If $Q$ has degree at
most $k$ on $S$, then $Q\circ\tau$ has degree at most $\ell k$ on $K$,
and since $f=g\circ\tau$, its error is at most $\|g-Q\|_{\infty,S}$. Taking infima proves \eqref{product_sandwich}.

The right inequality immediately transfers the rate from $f$ to $g$.
Conversely, if $E_k^S(g)=O(k^{-\beta})$, then for $k\geq2\ell$ put
$k'=\lfloor k/\ell\rfloor\geq k/(2\ell)$.  Monotonicity and the
left inequality give
\[
E_k^K(f)\leq E_{\ell k'}^K(f)\leq E_{k'}^S(g)=O(k^{-\beta}).
\]

For $k < 2 \ell$ the bound is trivial.

\end{proof}

\begin{proof}[Proof of Theorem~\ref{multi_dice}]
\emph{Necessity.} Suppose a multi-die enterprise for $f$ has the
uniform tail in (i). Run it using a single die of arbitrary weighting
$\mathbf z\in S$: whenever it requests die $i$, draw a fresh outcome of
the joint die and report only its $i$th coordinate, discarding the
others. At any history at which die $i$ is requested, the fresh joint
outcome is independent of the past, so the reported value is $a$ with
conditional probability $\mu(\mathbf z)_{i,a}$. This holds even though
the coordinates of a joint outcome may be correlated, and even though
the choice of $i$ is adaptive. By induction on the history, the
simulated run has the same law as a run of the original enterprise on
input $\mu(\mathbf z)$; in particular the output and the query count
have the same joint distribution. Fresh sampling is essential: reusing the unreported
coordinates of an earlier joint outcome could reveal correlations, and
would not simulate independent rolls.

This gives an enterprise for $g$ on all of $S$ that makes one joint
query per original query, so its query count has the same distribution and the
same uniform tail. Theorem~\ref{dice_approx_equiv} gives
$E_k^S(g)=O(k^{-2q})$, and Lemma~\ref{product_approx} gives
$E_k^K(f)=O(k^{-2q})$.

\emph{Sufficiency.} Suppose $E_k^K(f)=O(k^{-2q})$. By
Lemma~\ref{product_approx}, $E_k^S(g)=O(k^{-2q})$, so
Theorem~\ref{dice_approx_equiv} gives an oblivious enterprise for $g$
with uniform tail $O((1+n)^{-q})$. Answer each of its queries by rolling
each original die once, in the order $1,\ldots,\ell$, and returning the
tuple of outcomes. These tuples are independent with weighting
$\tau(\mathbf p_1,\ldots,\mathbf p_\ell)$, so the output is $1$ with
probability $g\circ\tau=f$. If $N_g$ is the budget chosen by the
enterprise for $g$, the simulation makes exactly $N_f=\ell N_g$ queries,
following a schedule of $N_g$ repetitions of $1,\ldots,\ell$ that is
fixed before any roll. Hence
\[
\PP(N_f>n)=\PP(N_g>\lfloor n/\ell\rfloor)
\leq C(1+\lfloor n/\ell\rfloor)^{-q}
\leq C\ell^q(1+n)^{-q},
\]
and $\E N_f=\ell\,\E N_g$ is uniformly bounded when $q>1$.
\end{proof}

This reduction controls the attainable exponent, not the dependence on
dimension: the joint die has $\prod_i(m_i+1)$ faces.

\section{Smoothness consequences}\label{sphere_interpretation}

Condition (ii) of Theorem~\ref{dice_approx_equiv} can be restated as a
smoothness condition after lifting $f$ to the sphere by the
substitution $p_i=u_i^2$. Let
$S^m=\{\mathbf u\in\mathbb R^{m+1}:\sum_i u_i^2=1\}$
and $\Phi(\mathbf u)=f(u_0^2,\ldots,u_m^2)$. For every continuous
$f:\simplex\to\mathbb R$,
\begin{equation}\label{square_identity}
E_d^{S^m}(\Phi)=E_{\lfloor d/2\rfloor}^{\simplex}(f)
\qquad(d\geq0).
\end{equation}
Indeed, if $P$ has degree at most $\lfloor d/2\rfloor$ on the simplex,
then $P(u_0^2,\ldots,u_m^2)$ has degree at most $d$ and, since
$\mathbf u\mapsto(u_0^2,\ldots,u_m^2)$ maps $S^m$ onto $\simplex$, the
same error. Conversely, average a degree-$d$ approximant of $\Phi$ over
the $2^{m+1}$ coordinate sign changes. Since $\Phi$ is invariant under
them, the error does not increase, and only monomials with all
exponents even survive; so the average has the form
$P(u_0^2,\ldots,u_m^2)$ with $\deg P\leq\lfloor d/2\rfloor$. (For
$m=0$, both sides of \eqref{square_identity} vanish.)

For $m\geq1$, let $R_{ij,\theta}$ be the rotation by angle $\theta$ in
the $(u_i,u_j)$-plane, put $T_{ij,\theta}G=G\circ R_{ij,\theta}$, and
define the modulus of smoothness
\[
\omega_r(G,t)_\infty
=\max_{0\leq i<j\leq m}\sup_{|\theta|\leq t}
\|(I-T_{ij,\theta})^rG\|_{\infty,S^m}.
\]

\begin{corollary}[Angular difference criterion]\label{smoothness}
Let $m\geq1$ and $q>0$, and let $r>2q$ be an integer. A continuous
$f:\simplex\to(0,1)$ has an enterprise with uniform tail
$O((1+n)^{-q})$ if and only if
\begin{equation}\label{angular_condition}
\omega_r(\Phi,t)_\infty=O(t^{2q})\qquad(t\downarrow0).
\end{equation}
\end{corollary}

\begin{proof}
Put $s=2q$. For $0<s<r$ and $G\in C(S^m)$, the direct and inverse
theorems of approximation on the sphere give
\[
E_k^{S^m}(G)=O(k^{-s})
\quad\Longleftrightarrow\quad
\omega_r(G,t)_\infty=O(t^s).
\]
For $m\geq2$, this is due to Dai and Xu
\cite[Theorem 3.4, (3.9)--(3.10); cf.\ Corollary 3.5]{dai2010moduli},
whose modulus in Definition 2.2 is the one above; their convention of
indexing degrees by $n-1$ does not affect rates. For $m=1$, spherical
polynomials are trigonometric polynomials and rotations are
translations, so this is the classical periodic equivalence
\cite[\S1, conditions (i$_1$)--(i$_2$)]{kolomoitsev2020smoothness}.
Apply the equivalence to $G=\Phi$, and then use \eqref{square_identity}
and Theorem~\ref{dice_approx_equiv}.
\end{proof}
Classical H\"older and Lipschitz regularity give a simple
sufficient condition, including at integer smoothness orders.

\begin{corollary}[A smooth sufficient condition]
\label{smooth_extension}
Let $K$ be a simplex or a finite product of simplices,
and let $s=a+\alpha$, where $a\geq0$ is an integer and
$0<\alpha\leq1$. Suppose $f:K\to(0,1)$ extends to a
$C^{a,\alpha}$ function on a neighborhood of $K$ in its
affine hull. Here $C^{a,1}$ means that the derivatives
of order $a$ are Lipschitz continuous.

Then $f$ has an oblivious enterprise with uniform query
tail $O((1+n)^{-s/2})$. If $s>2$, its expected query
count is uniformly bounded. In particular, a
$C^{2,\alpha}$ extension gives tail exponent
$1+\alpha/2$ for $0<\alpha\leq1$, and a $C^3$
extension gives exponent $3/2$.
\end{corollary}

\begin{proof}
A singleton domain needs no queries. Otherwise, $K$
is a convex polytope with nonempty interior in its
affine hull.

Consider any chord $I$ of $K$, parametrized affinely
by $[-1,1]$. Since the chord lengths are bounded by
$\operatorname{diam}(K)$, the $C^{a,\alpha}$ norms
of these parametrized restrictions of $f$ are bounded
uniformly over all chords. The univariate Jackson
inequality therefore gives
\[
E_k^I(f|_I)\lesssim k^{-s}
\]
with a constant independent of $I$. This applies
also when $\alpha=1$, since Lipschitz continuity
of the $a$th derivative gives the corresponding
integer-order approximation rate.

Totik's chordwise-to-polytope approximation result
\cite[\S1, paragraph preceding Remark~1.2]{totik2020polynomial} now gives
\[
E_k^K(f)=O(k^{-s}).
\]
Applying Theorem~\ref{dice_approx_equiv}, or
Theorem~\ref{multi_dice} for a product, with
$q=s/2$ proves the claim. A $C^3$ extension
has Lipschitz second derivatives on the compact
domain, so it satisfies the hypothesis with
$a=2$ and $\alpha=1$.
\end{proof}

Corollary~\ref{smooth_extension} is only a sufficient
condition. For a simplex, the exact criterion is
\eqref{angular_condition} for the spherical lift
$\Phi$. At integer $s=2q$, the exact condition is
of Zygmund type; the classical $C^{s-1,1}$
condition above is sufficient but not necessary.

For a single input coin and noninteger $s$, the
existence of a factory with uniform
$O((1+n)^{-s/2})$ query tail also follows from
Holtz et al.~\cite[Theorem~8]{holtz2011newcoins},
since $\Delta_n(p)\le n^{-1/2}$ uniformly in $p$.
The additional conclusions here are an oblivious
realization, the integer-order sufficient conditions,
and the extension to dice and products of simplices.

\section*{AI disclosure}
ChatGPT was used for mathematical exploration, proof development,
literature review, and manuscript editing. Claude and Codex also
assisted with manuscript editing. The author is responsible for the mathematical statements,
proofs, references, and presentation.

\appendix
\section{Translation of the weighted estimates}\label{coordinates}

Ditzian and Zhou~\cite{ditzian1993optimal} work in the coordinates
$x_i=p_i$, $1\leq i\leq m$, with $p_0=1-\sum_i x_i$. In these
coordinates, our weighted derivative along an edge $\{0,i\}$ is, up to
sign, $(x_ip_0)^{r/2}\partial_i^r$, and along an edge $\{i,j\}$ with
$i,j\geq1$ it is $(x_ix_j)^{r/2}(\partial_i-\partial_j)^r$. These are
exactly their weighted derivatives along unit edge directions: the
factors $\sqrt2$ in our directions and weights cancel. In particular,
the operator \eqref{diffusion_operator} becomes
\begin{equation}\label{coordinate_operator}
\mathcal LP
=\sum_{i=1}^m x_i(1-x_i)\partial_{ii}P
-2\sum_{1\leq i<j\leq m}x_ix_j\partial_{ij}P.
\end{equation}
The change of coordinates preserves total degree and uniform norm, so
the cited estimates apply on the full closed simplex.

\end{document}